\documentclass{amsart}

\usepackage{ae}

\usepackage[alphabetic]{amsrefs}
\usepackage{mathrsfs}

\renewcommand{\setminus}{\ \rule[2.5pt]{7pt}{1pt}\ }

\newcommand{\lie}[1]{\mathfrak{#1}}
\newcommand{\glie}{\lie{g}}
\newcommand{\hlie}{\lie{h}}

\newcommand{\plie}{\lie{p}}

\newcommand{\Lie}{\operatorname{Lie}}
\newcommand{\Stab}{\operatorname{Stab}}

\newcommand{\G}{\mathbf{G}}
\newcommand{\Gm}{{\G_m}}

\newcommand{\Aff}{\mathbf{A}}

\newcommand{\Ad}{\operatorname{Ad}}

\newcommand{\Int}{\operatorname{Int}}

\newcommand{\norm}[1]{\|#1\|}

\newcommand{\X}{\mathbf{X}}

\newcommand{\gcr}{$G$-cr}

\newcommand{\Z}{\mathbf{Z}}

\newcommand{\SL}{\operatorname{SL}}
\newcommand{\GL}{\operatorname{GL}}

\newtheorem{prop}{Proposition}
\newtheorem{lem}[prop]{Lemma}
\newtheorem{theorem}[prop]{Theorem}
\newtheorem{cor}[prop]{Corollary} 

\theoremstyle{remark}

\newtheorem{rem}[prop]{Remark}

\numberwithin{equation}{subsection}

\begin{document}
\title{Completely reducible Lie subalgebras}
\author{George McNinch}
\address{Department of Mathematics \\ 
         Tufts University \\ 
         503 Boston Avenue \\ 
         Medford, MA 02155 \\ 
         USA}
\email{george.mcninch@tufts.edu}
\date{December 30, 2005}
\thanks{Research of the author supported in part by the US National
  Science Foundation through DMS-0437482}

\begin{abstract}
  Let $G$ be a connected and reductive group over the algebraically
  closed field $K$. J-P. Serre has introduced the notion of a
  $G$-completely reducible subgroup $H \subset G$.  In this note, we
  give a notion of $G$-complete reducibility -- \gcr\ for short -- for
  Lie subalgebras of $\Lie(G)$, and we show that if the closed
  subgroup $H \subset G$ is \gcr, then $\Lie(H)$ is \gcr\ as well.
\end{abstract}

\maketitle

\section{Introduction}

Let $G$ be a connected and reductive group over the algebraically
field $K$, and write $\glie$ for the Lie algebra of $G$.  J-P. Serre
has introduced the notion of a $G$-completely reducible subgroup; we
state the definition here only for a closed subgroup $H \subset G$. We
say $H$ is \gcr\ provided that whenever $H \subset P$ for a parabolic
subgroup of $G$, there is a Levi factor $L \subset P$ such that $H
\subset L$; cf.  \cite{serre-sem-bourb}. When $G= \GL(V)$, the
subgroup $H$ is \gcr\ if and only if $V$ is a semisimple $H$-module.
Similarly, if the characteristic of $K$ is not 2 and $G$ is either the
symplectic group $\operatorname{Sp}(V)$ or the orthogonal group
$\operatorname{SO}(V)$, a subgroup $H$ of $G$ is \gcr\ if and only $V$
is a semisimple $H$-module.

B. Martin \cite{martin} used some techniques from ``geometric
invariant theory'' -- due to G.  Kempf and to G.  Rousseau -- to prove
that if $H \subset G$ is \gcr, and if $N$ is a normal subgroup of $H$,
then $N$ is \gcr\ as well; cf.  \cite{serre-sem-bourb}*{Th\'eor\`eme
  3.6}.  Martin's result was obtained first for \emph{strongly
  reductive} subgroups in the sense of Richardson; it follows from
\cite{bmr} that the strongly reductive subgroups of $G$ are precisely
the \gcr\ subgroups. See also \cite{serre-sem-bourb}*{\S 3.3} for an
overview of these matters.

We are going to prove in this note a result related to that of Martin.
If $\hlie \subset \glie$ is a Lie subalgebra, say that $\hlie$ is
\gcr\ provided that whenever $\hlie \subset \Lie(P)$ for a
parabolic subgroup $P$ of $G$, there is a Levi factor $L \subset
P$ such that $\hlie \subset \Lie(L)$.

We will prove:
\begin{theorem}
  \label{theorem:main-theorem}
  Suppose that $G$ is a reductive group over the algebraically closed
  field $K$.
  \begin{enumerate}
  \item Let $X_1,\dots,X_d$ be a basis for the Lie subalgebra $\hlie
    \subset \glie$. Then $\hlie$ is \gcr\ if and only if the $\Ad(G)$-orbit
    of $(X_1,\dots,X_d)$ is closed in $\bigoplus^d \glie$.
  \item If the closed subgroup $H \subset G$ is \gcr, then $\Lie(H)$ is \gcr\
    as well.
  \end{enumerate}
\end{theorem}

Our result -- and our techniques -- are related to those used by Richardson
in \cite{richardson}, though he treats mainly the case
of characteristic 0. See e.g. \emph{loc. cit.} Theorem 3.6.

The converse to Theorem \ref{theorem:main-theorem}(2) is not true.
Indeed, suppose the characteristic $p$ of $K$ is positive, and
consider a finite subgroup $H \subset G$ whose order is a power of
$p$. Then $\Lie(H) = 0$ is clearly \gcr; however, if $G=\SL(V)$ and if
$H$ is non-trivial, then $V$ is not semisimple as an $H$-module, thus
$H$ is not \gcr.  The converse to Theorem
\ref{theorem:main-theorem}(2) is even false for connected $H$; I thank
Ben Martin for pointing out the following example. Take for $H$ any
semisimple group in characteristic $p>0$, let $\rho_i:H \to \SL(V_i)$
be representations for $i=1,2$ with $\rho_1$ semisimple and $\rho_2$
not semisimple, and consider the representation $\rho:H \to G=\SL(V_1
\oplus V_2)$ given by $h \mapsto \rho_1(h) \oplus \rho_2(F(h))$ where
$F:H \to H$ is the Frobenius endomorphism. If $J$ denotes the image of
$\rho$, then $J$ is not \gcr\ since $V_1 \oplus V_2$ is not semisimple
as a $J$-module.  However, $\Lie(J) = \operatorname{im} d\rho$ lies in
the the Lie algebra of the subgroup $M = \SL(V_1) \times \SL(V_2)$;
moreover, $M$ is a Levi factor of a parabolic subgroup of $G$, and
$\Lie(J) = \operatorname{im}d\rho_1 \oplus 0 \subset \lie{sl}(V_1)
\oplus \lie{sl}(V_2) = \Lie(M)$.  Since the image of $\rho_1 \times
1:H \to M$ is $M$-cr (use \cite{bmr}*{Lemma 2.12(i)}), the main result
of this paper implies $\Lie(J)$ to be $M$-cr\footnote{This can be seen
  more easily: it is straightforward to check that a Lie subalgebra
  $\hlie \subset \lie{sl}(V)$ is $\SL(V)$-cr if and only if $V$ is a
  semisimple $\hlie$-module.}, hence Lemma \ref{lem:gcr-levi} below
shows that $\Lie(J)$ is \gcr\ as well.

The author would like to thank Michael Bate, Benjamin Martin, and
Gerhard R\"ohrle for some comments and conversations which were useful
in preparing this note. He would also like to acknowledge the
hospitality of the Centre Interfacultaire Bernoulli at the \'Ecole
Polytechnique F\'ed\'erale de Lausanne during a visit in June 2005
during which a need for the result of this note became clear to the
authors of \cite{mcninch-testerman}.

\section{Preliminaries}
\label{sec:prelim}

We work throughout in the \emph{geometric} setting; thus, $K$ is an
algebraically closed field. A variety will mean a separated and
reduced scheme of finite type over $K$. The group $G$ will be a
connected and reductive algebraic group (over $K$). A closed subgroup
$H \subset G$ is in particular a subvariety of $G$ and so $H$ is
necessarily \emph{reduced} -- e.g.  if $G$ acts on a variety $X$ and
if $x \in X$, then $\Stab_G(x)$ will mean the reduced subgroup
determined by the ``abstract group theoretic'' stabilizer (even if the
$G$-orbit of $x$ is not separable).

\subsection{Closed orbits}
\label{sub:closed-orbits}
Let $X$ be an affine $G$-variety, let $x \in X$ and choose a maximal
torus $S \subset \Stab_G(x)$ of the stabilizer in $G$ of $x$. Let $L =
C_G(S)$; thus $L$ is a Levi factor of a parabolic subgroup of $G$.

\begin{prop}
  \label{prop:closed-orbits}
  If the $G$-orbit $G \cdot x$ is closed in $X$, then the $L$-orbit $L
  \cdot x$ is closed in $X$.
\end{prop}

\begin{proof}
  The fixed point set $X^S$ is closed in $X$. Since by assumption $G
  \cdot x$ is closed in $X$, it follows that
  \begin{equation*}
    (G \cdot x)^S = X^S \cap G \cdot x
  \end{equation*}
  is closed in $X$.

  Let now $N = N_G(S)$ be the normalizer in $G$ of $S$.
  We claim that $(G \cdot x)^S = N \cdot x$. Indeed, let $g \in G$
  and suppose $g \cdot x$ is fixed by $S$. The claim follows once we
  prove that $g \cdot x \in N \cdot x$. Well, for each $s \in S$ we
  have $sg \cdot x = g \cdot x$ so that $g^{-1}sg \in \Stab_G(x)$.
  Thus $g^{-1}Sg$ is a maximal torus of $\Stab_G(x)$. Since maximal
  tori are conjugate \cite{springer-LAG}*{Theorem 6.4.1}, there is an
  element $h \in\Stab_G(x)$ such that $g^{-1}Sg = hSh^{-1}$. But then
  $gh \in N$, and moreover, $g \cdot x = gh \cdot x$. 

  Note that $N$ contains $L$ as a normal subgroup.  We now observe
  that the stabilizer in $L$ of a point $y$ of the orbit $N \cdot x$
  is conjugate to $\Stab_L(x)$ by an element of $N$. Indeed, choosing
  $h \in N$ such that $h \cdot x = y$, one knows that $h \cdot
  \Stab_L(y) \cdot h^{-1} = \Stab_L(x)$.  It follows that all
  $L$-orbits in $N \cdot x$ have the same dimension.

  Since the closure of any $L$-orbit must be the union of orbits of
  strictly smaller dimension, it follows that the $L$-orbits in $N
  \cdot x$ are closed.

  Since $N \cdot x = X^S \cap G \cdot x$ is closed in $X$, it follows
  at once that $L \cdot x$ is closed in $X$, as required.
\end{proof}

\begin{rem}
  With notation as in the previous Proposition, if $N = N_G(S)$, it
  follows from the rigidity of tori \cite{springer-LAG}*{Corollary
    3.2.9} that $L$ has finite index in $N$.  In particular, $(G \cdot
  x)^S$ is a finite union of $L$-orbits which are permuted
  transitively by $N$; moreover, these $L$-orbits are precisely the
  connected components of $(G \cdot x)^S$.
\end{rem}

\subsection{Complete reducibility}
The interpretation of complete reducibility using the spherical
building of $G$ permits one to prove the following:
 \begin{lem}
   \label{lem:gcr-levi}
   Let $G$ be reductive and let $M \subset G$ be a Levi factor of a
   parabolic subgroup of $G$. Suppose that $J \subset M$ is a
   subgroup, and that $\hlie \subset \Lie(M)$ is a Lie subalgebra.
   Then $J$ is \gcr\ if and only if $J$ is $M$-cr and $\hlie$ is \gcr\
   if and only if $\hlie$ is $M$-cr.
 \end{lem}

 \begin{proof}
   The assertion for $J$ follows from
   \cite{serre-sem-bourb}*{Proposition 3.2}. The proof for $\hlie$ is
   similar; let us give a sketch. Write $X$ for the building of
   $G$.  The Lie subalgebra $\hlie$ defines a subcomplex $Y$ of
   $X$: the simplices of $Y$ are those simplices in $X$ which
   correspond to parabolic subgroups $P$ with $\hlie \subset \Lie(P)$.

   Recall \cite{borel-LAG}*{Corollary 14.13} that the intersection $P
   \cap P'$ of two parabolic subgroups $P,P' \subset G$ contains a
   maximal torus of $G$.  This implies that $\Lie(P \cap P') = \Lie(P)
   \cap \Lie(P')$; see e.g. the argument in the first paragraph of
   \cite{jantzen-nilpotent}*{\S 10.3}.

   If now $\hlie \subset \Lie(P) \cap \Lie(P')$, it follows that
   $\hlie \subset \Lie(P \cap P')$. This shows that the subcomplex $Y$
   is convex; see \cite{serre-sem-bourb}*{Prop. 3.1}.  Evidently
   $\hlie$ is \gcr\ if and only if $Y$ is $X$-cr in the sense of
   \cite{serre-sem-bourb}*{\S2.2}.

   Choose a parabolic subgroup $Q$ for which $M$ is a Levi factor.
   Then we may identify the building of $M$ with the residual
   building of $X$ determined by the parabolic $Q$; cf.
   \cite{serre-sem-bourb}*{2.1.8 and 3.1.7}.  Now the claim follows from
   \cite{serre-sem-bourb}*{Proposition 2.5}.
 \end{proof}

\subsection{Cocharacters and parabolic subgroups}

If $V$ is a variety and $f:\G_m \to V$ is a morphism, we write
$v = \lim_{t\to 0} f(t)$, and we say that the limit
exists, if $f$ extends to a morphism $\tilde f:\Aff^1 \to V$ with
$\tilde f(0) = v$.  If $\Gm$ acts on $V$, a closed point $w \in V$
determines a morphism $f:\Gm \to V$ via the rule $t \mapsto t\cdot w$;
one writes $\lim_{t \to 0} t \cdot w$ as shorthand for $\lim_{t\to 0}
f(t)$.

A \emph{cocharacter} of an algebraic group $A$ is a $K$-homomorphism
$\gamma:\Gm \to A$. A linear $K$-representation $(\rho,V)$ of $A$
yields a linear $K$-representation $(\rho \circ \gamma,V)$ of $\Gm$.
Then $V$ is the direct sum of the weight spaces
\begin{equation}
  \label{eq:weight-spaces}
  V(\gamma;i) = \{v \in V \mid (\rho \circ \gamma)(t)v = t^i v, \ \forall\ t \in \Gm \}
\end{equation}
for $i \in \Z$. We write $X_*(A)$ for the set of cocharacters of $A$.

Consider now the  reductive group $G$.  If $\gamma \in X_*(G)$,
then
\begin{equation*}
  P_G(\gamma)= P(\gamma) = 
  \{ x \in G \mid \lim_{t \to 0} \gamma(t)x\gamma(t^{-1})\text{\ exists}\}
\end{equation*}
is a parabolic subgroup of $G$ whose Lie algebra is $\plie(\gamma) =
\sum_{i \ge 0} \glie(\gamma;i)$.  Moreover, each parabolic subgroup of
$G$ has the form $P(\gamma)$ for some cocharacter $\gamma$; for all
this cf. \cite{springer-LAG}*{3.2.15 and 8.4.5}.  

We note that $\gamma$ ``exhibits'' a Levi decomposition of
$P=P(\gamma)$. Indeed, $P(\gamma)$ is the semi-direct product
$Z(\gamma) \cdot U(\gamma)$, where $U(\gamma) = \{x \in P \mid
\lim_{t\to 0} \gamma(t)x\gamma(t^{-1}) = 1\}$ is the unipotent radical of
$P(\gamma)$, and the reductive subgroup $Z(\gamma) =
C_G(\gamma(\G_m))$ is a Levi factor in $P(\gamma)$; cf.
\cite{springer-LAG}*{13.4.2}.

\begin{lem}
  \label{lem:levi-cocharacter}
  Let $P$ be a parabolic subgroup of $G$, let $L$ be a Levi factor of
  $P$, let $\gamma \in X_*(L)$ and assume that $P=P(\gamma)$.  Then
  $L=Z(\gamma)$ and the image of $\gamma$ lies in the connected center
  of $L$.
\end{lem}

\begin{proof}
  Let $R$ be the radical of $P$. Then the Levi factors of $P$ are
  precisely the centralizers of the maximal tori of $R$; cf.
  \cite{borel-LAG}*{Cor. 14.19}.  Since the connected center of a Levi
  factor of $P$ evidently lies in $R$, we see that the connected
  center of each Levi factor is a maximal torus of $R$.

  Now, the centralizer $L_1=Z(\gamma)$ is a Levi factor of $P$, so
  that $\gamma$ is a cocharacter of the connected center of $L_1$; in
  particular, the image of $\gamma$ lies in $R$.  Moreover, since
  $L_1=Z(\gamma)$, the centralizer of the image of $\gamma$ in $R$ is
  a maximal torus $S$ of $R$. It follows that $S$ is the unique
  maximal torus of $R$ containing the image of $\gamma$.

  Since the image of $\gamma$ lies in $L$ and in $R$, and since $L$
  intersects $R$ in a maximal torus of $R$, it follows that $S = L
  \cap R$ so that $L = L_1$ as required.
\end{proof}

\subsection{Instability in invariant theory}
Let $(\rho,V)$ be a linear representation [always assumed finite
dimensional] of $G$, and fix a closed $G$-invariant subvariety $S
\subset V$.  We are going to describe a precise form -- due to Kempf
and Rousseau -- of the Hilbert-Mumford criteria for the instability of
a vector $v \in V$ under the action of $G$.

Let us first briefly describe our goal: given a Lie subalgebra $\hlie
\subset \glie = \Lie(G)$, fix a basis $\mathbf{X} = (X_1,\dots,X_d)$
of $\hlie$. If the $G$-orbit of $\mathbf{X}$ in $\bigoplus^d \glie$ is
not closed -- so that $\mathbf{X}$ is an unstable vector -- the
results of Kempf and Rousseau permit us to associate to $\mathbf{X}$ a
unique parabolic subgroup $P_{\mathbf{X}}$; see Corollary
\ref{cor:not-closed-case} below. If $g \in G$ satisfies $\Ad(g)\hlie =
\hlie$, one of our main objectives is to show that $g \in
P_{\mathbf{X}}$. Using $g$, we get a new basis $\Ad(g)\mathbf{X} =
(\Ad(g)X_1,\dots,\Ad(g)X_d)$ of $\hlie$, and generalities show that
$P_{\Ad(g)\mathbf{X}} = gP_{\mathbf{X}}g^{-1}$. So we want to prove
the equality $P_{\mathbf{X}} = P_{\Ad(g)\mathbf{X}}$; it will then
follow that $g\in P_{\mathbf{X}}$, as desired.

Return now to our general setting: $V$ is any linear representation of
$G$.  For $v \in V$, put
\begin{equation*}
 |V,v| = \{ \lambda \in X_*(G) \mid  \lim_{t \to 0}
 \rho(\lambda(t))v \quad \text{exists} \}.
\end{equation*}
Write $V = \bigoplus_{i \in \Z} V(\lambda;i)$ as in \eqref{eq:weight-spaces},
and write $v = \sum_i v_i$ with $v_i \in V(\lambda;i)$.
Then evidently
\begin{equation}
  \label{eq:positive-cond}
  \lambda \in |V,v| \iff v_i = 0 \quad \forall i < 0;
\end{equation}
if $\lambda \in |V,v|$ then of course $\lim_{t \to 0}
\rho(\lambda(t))v = v_0$.

Now let $S \subset V$ be a $G$-invariant closed subvariety and suppose
that $v \not \in S$.  Given $\lambda \in |V,v|$, write $v_0 =
\lim_{t\to 0} \rho(\lambda(t))v$. If $v_0 \in S$, write
$\alpha_{S,v}(\lambda)$ for the order of vanishing of the regular
function $(t \mapsto \rho(\lambda(t))v - v_0):\Aff^1 \to V$, otherwise
write $\alpha_{S,v}(\lambda) = 0$; see \cite{kempf-instab}*{\S3} for
more details.  Then $\alpha_{S,v}(\lambda)$ is a non-negative integer,
and $\alpha_{S,v}(\lambda) > 0$ if and only if $v_0 \in S$.  Moreover,
if $v=\sum_{i \in \Z} v_i$ with $v_i \in V(\lambda;i)$ as before, then
\begin{equation}
  \label{eq:alpha-descr}
  v_0 \in S \ \implies \  \alpha_{S,v}(\lambda) = 
  \alpha_{\{v_0\},v}(\lambda) =  \min\{j>0 \mid v_j \ne 0\}.
\end{equation}

Suppose that $W \subset V$ is a subspace of dimension $d = \dim W$.
Let $w_1,\dots,w_d$ be a basis of $W$, and consider the point $x =
(w_1,\dots,w_d)$ of the linear space $X = \bigoplus^d V$; abusing
notation somewhat, we write also $\rho$ for the diagonal action
$\bigoplus^d \rho$ of $G$ on $X$. We observe  for $\lambda \in
X_*(G)$ that we have
\begin{equation}
  \label{eq:basis-equiv-cond}
  \lambda \in |X,x| \iff W \subset \sum_{j \ge 0} V(\lambda;j).
\end{equation}

% \begin{lem}
%   \label{lem:positive-lemma}
%   \begin{enumerate}
%   \item $\lambda \in |V,v|$ if and only if $v_i =0$ for each $i<0$.
%   \item If $\lambda \in |V,v|$, then $\lim_{t\to 0} \rho(\lambda(t))v
%     = v_0$. 
%   \item If $v_0 \in S$, then
%     \begin{equation*}
%       \alpha_{S,v}(\lambda)
%       = \alpha_{\{v_0\},v}(\lambda)
%       = \min\{j>0 \mid v_j \ne 0\};
%     \end{equation*}
%     otherwise, $\alpha_{S,v}(\lambda)=0$.
%   \end{enumerate}
% \end{lem}

% \begin{proof}
%   Immediate from the definitions.
% \end{proof}

% \begin{lem}
%   \label{lem:basis-equiv-cond}
%   Let $w_1,\dots,w_d$ be a basis for $W$, and let $x=(w_1,\dots,w_d)
%   \in X=\bigoplus^dV$. 
%   The following are equivalent:
%   \begin{enumerate}
%   \item $\lambda \in |X,x|$.
%   \item $W \subset \sum_{j \ge 0} V(\lambda;j)$.
%   \item $\lim_{t \to 0} \rho(\lambda(t))v$ exists for each $v \in W$.
%   \end{enumerate}
% \end{lem}

% \begin{proof}
%   Since $\lim_{t \to 0} \rho(\lambda(t))x$ exists if and only if
%   $\lim_{t \to 0} \rho(\lambda(t))w_i$ exists for $1 \le i \le d$, the
%   equivalence of these three conditions is a consequence of Lemma
%   \ref{lem:positive-lemma}.
% \end{proof}

Fix $S \subset X=\bigoplus^d V$ a closed and $\rho(G)$-invariant
subvariety, and assume that $x = (w_1,\dots,w_d) \not \in S$. In
this setting one may compute the function $\alpha_{S,x}$ for the
diagonal $G$-action on $X$ using functions $\alpha_{\{v_0\},v}$ for
the $G$-representation $V$. More precisely, we have:
\begin{lem}
  \label{lem:alpha-indep}
  Let $\lambda \in |X,x|$ and suppose 
  $\alpha_{S,x}(\lambda)>0$. For $w \in W$, write $w_0 = \lim_{t \to
    0} \rho(\lambda(t))w$. Then
  \begin{equation*}
    (*) \quad \alpha_{S,x}(\lambda) = \min_{w \in W} \alpha_{\{w_0\},w}(\lambda).
  \end{equation*}
\end{lem}

\begin{proof}
  For $1 \le i \le d$ write $x = \sum_j x^j$ with $x^j \in
  X(\lambda;j)$. 

% Moreover, for $j \in \Z$, write
%   \begin{equation*}
%     x^j = (w_1^j,\dots,w_d^j) \quad \text{for}\quad w_i^j \in V(\lambda;j), \quad 1 \le i \le d.
%   \end{equation*}
%   Thus $w_i = \sum_j w_i^j$ for $1 \le i \le d$.

  By assumption, $\lambda \in |X,x|$; by 
  \eqref{eq:positive-cond}  we see that $x^j = 0$ if $j < 0$.
  Moreover, using \eqref{eq:alpha-descr} we see that
  \begin{equation}
    \label{eq:alpha-1}
    \alpha_{S,x}(\lambda) = \alpha_{\{x^0\},x}(\lambda) = \min(j>0 \mid x^j \ne 0).
  \end{equation}
  If we now write $\displaystyle R = \min_{v \in W}
  \alpha_{\{v_0\},v}(\lambda)$ for the right hand side of $(*)$, then
  upon considering the components in $V$ of the vectors $x^j \in X =
  \bigoplus^d V$, one uses \eqref{eq:alpha-1} to see that
  $\alpha_{S,x}(\lambda) \ge R$.

  On the other hand, we may choose $v \in W$ such that $R =
  \alpha_{\{v_0\},v}(\lambda).$ Writing $v = \sum_{j \ge 0} v^j$ with
  $v^j \in V(\lambda;j)$, we see that
  \begin{equation*}
    \label{eq:R-eq}
    R=    \alpha_{\{v^0\},v}(\lambda) = \min(j > 0 \mid v^j \ne 0)
  \end{equation*}
  by \eqref{eq:positive-cond}.  Now write
  \begin{equation*}
    v = \sum_i \beta_i w_i \quad \text{  for scalars $\beta_i \in K$.}
  \end{equation*}
  Now, $v^R \ne 0$ implies that $x^R \ne 0$; it
  follows from \eqref{eq:alpha-1} that $R \ge
   \alpha_{S,x}(\lambda)$, and the Lemma is proved.
\end{proof}

Fix a basis $\{w_i\}$ for $W$ and let $x = (w_1,\dots,w_d) \in X$. Write
\begin{equation*}
  S = \overline{\rho(G)x} \setminus \rho(G)x;
\end{equation*}
then $S$ is closed in $X$ Notice that $S$ is a closed subset, since
$\rho(G)x$ is open in $\overline{\rho(G)x}$, and $S$ is $G$-invariant.
We suppose that $\rho(G)x$ is not closed, or equivalently that $S$ is
non-empty.

\begin{cor}
  \label{cor:different-bases}
  Let $h \in G$ satisfy $\rho(h)W = W$. If $x' = \rho(h)x$, then we
  have $|X,x| = |X,x'|$.  Moreover,
  \begin{equation*}
    \alpha_{S,x}(\lambda) = \alpha_{S,x'}(\lambda)
  \end{equation*}
  for each $\lambda \in |X,x|$.
\end{cor}

\begin{proof}
  Since by \eqref{eq:basis-equiv-cond} the sets $|X,x|$ and $|X,x'|$
  both consist of all cocharacters $\lambda$ for which $W \subset
  \sum_{j \ge 0} V(\lambda;j)$, we have that $|X,x| = |X,x'|$.

  Now write $x_0 = \lim_{t \to 0} \rho(\lambda(t))x$ and $x_0' =
  \lim_{t \to 0} \rho(\lambda(t))x'$.  We first claim that $x_0 \in S$
  if and only if $x_0' \in S$.

  Well, assume that $x_0 \not \in S$. Since $x_0$ lies in the closure
  of $\rho(G)x$ but not in $S$, it actually lies in $\rho(G)x$; thus
  $(\dagger) \ x_0 = \rho(g)x$ for some $g \in G$.

  Since the components in $V$ of the vector $x \in X = \bigoplus^d V$
  form a basis of $W$, one concludes from $(\dagger)$ that
  \begin{equation*}
    \lim_{t  \to 0}\rho(\lambda(t)) y = \rho(g)y
  \end{equation*}
  for each $y \in \bigoplus^d W \subset X$.  This shows in particular
  that $x_0' = \rho(g)x' = \rho(gh)x$, so that $x_0' \not \in S$. Since the
  argument just given is symmetric in $x$ and $x'$, it follows that
  $x_0 \in S$ if and only if $x_0' \in S$.

  Recall that $\alpha_{S,x}(\lambda) > 0$ if and only if $x_0 \in S$
  and that $\alpha_{S,x'}(\lambda)>0$ if and only if $x_0' \in S$. Thus to
  prove the final equality asserted by the corollary, we may suppose
  that $x_0,x_0' \in S$. Now, according to $(*)$ of Lemma
  \ref{lem:alpha-indep} we have
  \begin{equation*}
    \alpha_{S,x}(\lambda) = \min_{w \in W} \alpha_{\{w_0\},w}(\lambda) 
    = \alpha_{S,x'}(\lambda)
  \end{equation*}
  as required.
\end{proof}

Fix a real-valued $G$-invariant length function $\lambda \mapsto
\norm{\lambda}$ on the set $X_*(G)$ of cocharacters of $G$.

\begin{theorem}[Kempf \cite{kempf-instab}*{Theorem 3.4}, Rousseau]
  \label{theorem:kempf-rousseau}
  Let $z \in X \setminus S$ and assume that $\overline{\rho(G)z} \cap S$ is
  non-empty. Then the function $\alpha_{S,z}(\lambda)/\norm{\lambda}$
  assumes a maximal value $B>0$ on the non-trivial elements of
  $|X,z|$.  Let
  \begin{equation*}
    \Delta_{S,z} = \{\lambda \in |X,z| \mid
    \alpha_{S,z}(\lambda) = B \cdot \norm{\lambda} \quad \text{and
      $\lambda$ is indivisible}\}.
  \end{equation*}
  Then
  \begin{enumerate}
  \item  $\Delta_{S,z}$ is non-empty,
  \item there is a parabolic subgroup $P_{S,z}$ of $G$ such that
    $P_{S,z} = P(\lambda)$ for each $\lambda \in \Delta_{S,z}$,
  \item $\Delta_{S,z}$ is a principal homogeneous space under $R_uP_{S,z}$, and
  \item any maximal torus of $P_{S,z}$ contains a unique
    cocharacter which lies in $\Delta_{S,z}$.
  \end{enumerate}
\end{theorem}

Let $H \subset G$ be a  subgroup and suppose that $W$ is $\rho(H)$
invariant. Let $x = (w_1,\dots,w_d) \in X$ for a basis $\{w_i\}$ of $W$.
\begin{cor}
  \label{cor:not-closed-case}
  Assume that $\rho(G)x$ is not closed in $X$, and let 
  \begin{equation*}
    S = \overline{\rho(G)x} \setminus \rho(G)x. 
  \end{equation*}
  Then
  \begin{enumerate}
  \item $P_{S,x}$ is a proper parabolic subgroup of $G$, 
  \item $H \subset P_{S,x}$, and
  \item if $L \subset P_{S,x}$ is a Levi factor, there is a
    cocharacter $\lambda$ of the connected center $Z$ of $L$ which
    lies in $\Delta(S,x)$.
  \end{enumerate}
\end{cor}

\begin{proof}
  Since the image of any $\lambda \in |X,x|$ with $\alpha_{S,x}(\lambda) > 0$ is
  not central in $G$, (1) is immediate.

  Since the parabolic subgroup $P = P_{S,x}$ is self-normalizing, (2)
  will follow if we show that $hPh^{-1} = P$ for each $h \in H(k)$.
  But $hP_{S,x}h^{-1} = P_{S,\rho(h)x}$; see e.g.
  \cite{kempf-instab}*{Cor. 3.5}. Since $\rho(h)W = W$, Corollary
  \ref{cor:different-bases} shows that $|X,x| = |X,\rho(h)x|$ and that
  $\alpha_{S,x}(\lambda) = \alpha_{S,\rho(h)x}(\lambda)$ for all
  $\lambda \in |X,x| = |X,\rho(h)x|$; thus $\Delta_{S,x} =
  \Delta_{S,\rho(h)x}$ so that $P_{S,x} = P_{S,\rho(h)x}$ by Theorem
  \ref{theorem:kempf-rousseau}. Thus indeed $H \subset P_{S,x}$.

  Finally, for (3) let $S$ be a maximal torus of $L$ and hence of
  $P_{S,x}$. By (3) of Theorem \ref{theorem:kempf-rousseau}, $S$ has a
  cocharacter $\lambda$ which lies in $\Delta_{S,x}$.  Since $P_{S,x}
  = P(\lambda)$, it follows from Lemma \ref{lem:levi-cocharacter} that
  the image of $\lambda$ lies in the connected center of $L$, as
  required.
\end{proof}

Finally, we record:
\begin{lem}
  \label{lem:closed-case}
  Assume that $\rho(G)x$ is closed in $X$ and that $\lambda \in
  |X,x|$.  Then the subset 
  \begin{equation*}
    \lim_{t \to 0} \rho(\lambda(t))W =  
    \left \{\lim_{t \to 0} \rho(\lambda(t))w \mid w \in W \right\}
  \end{equation*}
  satisfies 
  \begin{equation*}
    \lim_{t \to 0} \rho(\lambda(t))W = \rho(g)W
  \end{equation*}
  for some $g \in G$.
\end{lem}

\begin{proof}
  Since $\lambda \in |X,x|$, the limit $x_\lambda = \lim_{t \to 0}
  \rho(\lambda(t))x$ exists. Since the orbit $\rho(G)x$ is closed, we
  have $\rho(g)x = x_\lambda$ for some $g \in G$.  Since
  $w_1,\dots,w_d$ is a basis of $W$, it follows that $\Ad(g)w =
  \lim_{t \to 0} \rho(\lambda(t))w$ for each $w \in W$, whence the
  Lemma.
\end{proof}

\section{Proof of the main theorem}

Recall that $G$ is a reductive group with Lie algebra $\glie$, and
that $\hlie \subset \glie$ is a Lie subalgebra.  Fix a basis
$X_1,\dots,X_d \in \hlie$, and let $\X = (X_1,\dots,X_d) \in \bigoplus^d
\glie = Y$.  We write $(\Ad,Y)$ for the representation
$(\bigoplus^d \Ad,\bigoplus^d \glie)$ of $G$.

\begin{proof}[Proof of part (1) of Theorem \ref{theorem:main-theorem}]
  Recall that we must show: the Lie algebra $\hlie$ is \gcr\ if and only if the
  $G$-orbit of $\X$ is closed in $Y=\bigoplus^d \glie$. 

  We first suppose that $\Ad(G)\X$ is closed, and we show that $\hlie$
  is \gcr.  Let $S$ be a maximal torus of the centralizer
  $C_G(\hlie)$. Then $\hlie \subset \Lie(L)$ where $L=C_G(S)$;
  moreover, $L$ is a Levi factor of a parabolic subgroup of $G$.  It
  follows from Lemma \ref{lem:gcr-levi} that $\hlie$ is \gcr\ if and
  only if $\hlie$ is $L$-cr.

  Moreover, it follows from Proposition \ref{prop:closed-orbits} that
  $\Ad(L)\X$ is closed in $Y$.  Thus we may replace $G$ by $L$ and so
  suppose that any torus in $G$ which centralizes $\hlie$ is central
  in $G$.  [Equivalently: $\hlie$ is not contained in the Lie algebra
  of any Levi factor of a proper parabolic subgroup of $G$.]  To show
  that $\hlie$ is \gcr\, we will show that $\hlie$ is not contained in
  $\Lie(P)$ for any proper parabolic subgroup $P$ of $G$. 

  Suppose that $\hlie \subset \Lie(P)$ for a parabolic subgroup
  $P\subset G$; we will show that $P=G$. Write $P = P(\phi)$ for some
  cocharacter $\phi$ of $G$, and write $L = L(\phi)$ for the
  centralizer in $G$ of the image of $\phi$; then $L$ is a Levi factor
  of $P$. 

  Since the $G$-orbit of $\mathbf{X}$ is closed, Lemma
  \ref{lem:closed-case} shows that 
  \begin{equation*}
    \lim_{t \to 0}    \Ad(\phi(t))\hlie = \Ad(g) \hlie 
  \end{equation*}
  for some $g \in G$. Since $\lim_{t\to 0} \Ad(\phi(t))H \in \Lie(L)$
  for each $H \in \hlie$, we conclude that $\hlie \subset
  \Ad(g^{-1})\Lie(L)$. But then the image of the cocharacter
  $\Int(g^{-1}) \circ \phi$ is a torus centralizing $\hlie$; hence the
  image of $\phi$ is central in $G$ so that $P=G$. This proves that
  $\hlie$ is indeed \gcr.

  To complete the proof of (i), it remains to show: if the orbit
  $\Ad(G)\X$ is not closed, then $\hlie$ is not \gcr. As
  in Corollary \ref{cor:not-closed-case} let $S = \overline{\Ad(G)\X}
  \setminus \Ad(G)\X$; our assumption means that $S$ is non-empty so
  that $\alpha_{S,\X}(\lambda) > 0$ for each $\lambda \in
  \Delta_{S,\X}$.  Moreover, $P = P_{S,\X}$ is a proper parabolic
  subgroup of $G$.

  We have $\hlie \subset \Lie(P)$ by 
  \eqref{eq:basis-equiv-cond}. To complete the proof, we suppose
  $\hlie$ is \gcr, and find a contradiction. 

  Since $\hlie$ is \gcr, there is a Levi factor $L$ of $P$ with $\hlie
  \subset \Lie(L)$. By Corollary \ref{cor:not-closed-case}, there is a
  cocharacter $\lambda$ of the connected center of $L$ which lies in
  $\Delta_{S,\X}$. Since $\hlie \subset \Lie(L)$, we have $\hlie
  \subset \glie(\lambda;0)$; thus $\X \in X(\lambda;0)$. But then
  $\alpha_{S,\X}(\lambda)=0$,  which is impossible since
  $\lambda \in \Delta_{S,\X}$.
\end{proof}

\begin{proof}[Proof of part (2) of Theorem \ref{theorem:main-theorem}]
  Recall that if $H \subset G$ is a subgroup which is \gcr, we
  must prove that $\hlie = \Lie(H)$ is \gcr.

  Let $S \subset C_G(H)$ be a maximal torus. Then $H \subset L =
  C_G(S)$ and $\hlie \subset \Lie(L)$. Applying
  Lemma \ref{lem:gcr-levi}, it is enough to show that
  $\hlie$ is $L$-cr; thus we replace $G$ by $L$ and so suppose that
  $H$ is not contained in a Levi factor of any proper parabolic
  subgroup of $G$. Since $H$ is \gcr, we conclude that $H$ is
  contained in no proper parabolic subgroup of $G$.

  To show that $\hlie$ is \gcr, we use part (1) of Theorem
  \ref{theorem:main-theorem}; it is enough to show that $\Ad(G)\X$ is
  closed in $Y$.  In fact, we are going to suppose that $\Ad(G)\X$ is
  not closed and obtain a contradiction.  Let $S=\overline{\Ad(G)\X}
  \setminus \Ad(G)\X$ and let $P=P_{S,\X}$.  Since $S$ is assumed
  non-empty, Corollary \ref{cor:not-closed-case} shows that $P$ is a
  proper parabolic subgroup. Moreover, since $\Ad(H)$ leaves $\hlie$
  invariant, that same corollary shows that $H \subset P$. This
  contradiction completes the proof.
\end{proof}

\newcommand\mylabel[1]{#1\hfil}

\end{document}